\documentclass[12pt,a4paper,reqno]{amsart}
\usepackage{amssymb}
\usepackage{latexsym}
\usepackage{exscale}
\usepackage{mathrsfs}

\usepackage{color}

\newcommand{\cH}{\mathcal H}
\newcommand{\cX}{\mathcal X}

\newcommand{\sB}{\mathscr B}

\newcommand{\be}{{\mathbf e}}

\newcommand {\SC} {{\mathbb C}}

\newcommand {\SH} {{\mathbb H}}
\newcommand {\SN} {{\mathbb N}}
\newcommand {\SR} {{\mathbb R}}
\newcommand {\ST} {{\mathbb T}}

\newcommand {\SX} {{\mathbb X}}

\newcommand {\SZ} {{\mathbb Z}}

\newcommand {\al} {{\alpha}}
\newcommand {\dt} {{\delta}}

\newcommand {\e} {{\varepsilon}}
\newcommand {\ga} {{\gamma}}

\newcommand {\la} {{\lambda}}
\newcommand {\La} {{\Lambda}}

\newcommand {\om} {{\omega}}

\newcommand{\fX}{{\mathfrak{X}}}

\def\supp{\mathop{\rm supp}}

\numberwithin{equation}{section}
\newtheorem{theorem}{Theorem}[section]
\newtheorem{lemma}[theorem]{Lemma}

\newtheorem{corollary}[theorem]{Corollary}
\newtheorem{Remark}[theorem]{Remark}
\newtheorem{proposition}[theorem]{Proposition}
\newtheorem{definition}[theorem]{Definition}
\newtheorem{example}[theorem]{Example}

\newcommand {\Proofof}[1] {\noindent{\bf P{\footnotesize\bf ROOF} of {#1}: } \ }

\newcommand {\Proof} {\noindent{\bf P{\footnotesize\bf ROOF}: } \ }
\newcommand {\ProofEnd} {
             \begin{flushright} \vskip -0.2in $\Box$ \end{flushright}}

\newcommand{\Ba}[1]{\begin{array}{#1}}
\newcommand{\Ea}{\end{array}}
\newcommand{\Be}{\begin{equation}}
\newcommand{\Ee}{\end{equation}}
\newcommand{\Bea}{\begin{eqnarray}}
\newcommand{\Eea}{\end{eqnarray}}
\newcommand{\Beas}{\begin{eqnarray*}}
\newcommand{\Eeas}{\end{eqnarray*}}
\newcommand{\Benu}{\begin{enumerate}}
\newcommand{\Eenu}{\end{enumerate}}
\newcommand{\Bi}{\begin{itemize}}
\newcommand{\Ei}{\end{itemize}}

\newcommand{\BR}{\begin{Remark} \em}
\newcommand{\ER}{\end{Remark}}
\newcommand{\BE}{\begin{example} \em}
\newcommand{\EE}{\end{example}}

\newcommand {\mand} {{\quad\mbox{and}\quad}}
\renewcommand {\mid} {{\,\,\,\colon\,\,\,}}

\newcommand{\bline}{{\bigskip

\noindent}}

\newcommand {\bone} {{\bf 1}}

\renewcommand {\span} {\mbox{\rm span}\,}

\newcommand{\xy}{{\langle x,y\rangle}}

\newcounter{reg}
\begin{document}

\title{Conditional quasi-greedy bases in Hilbert and Banach
spaces}

\author{G. Garrig\'os}
\address{Gustavo Garrig\'os
\\
Departamento de Matem\'aticas
\\
Universidad de Murcia
\\
30100 Murcia, Spain} \email{gustavo.garrigos@um.es}

\author{P. Wojtaszczyk}
\address{P. Wojtaszczyk, Interdisiplinary Centre for Mathematical and Computational Modelling,
University of Warsaw, 02-838 Warszawa,
ul. Prosta 69, Poland,
and
Institut of Mathematics, Polish Academy of Sciences\\00-956 Warszawa, ul. \`Sniadeckich 8, Poland.} \email{wojtaszczyk@mimuw.edu.pl}

\thanks{First author partially supported by grants MTM2010-16518 and MTM2011-25377 (Spain).
The second author was partially supported by the ``HPC
Infrastructure for Grand Challenges of Science and Engineering"
Project, co-financed by the European Regional Development Fund
under the Innovative Economy Operational Programme" and Polish NCN
grant DEC2011/03/B/ST1/04902.}



\maketitle

\begin{abstract}
For quasi-greedy bases $\sB$ in Hilbert spaces, we give
an improved bound of the associated conditionality constants
$k_N(\sB)=O(\log N)^{1-\e}$, for some $\e>0$, answering a question by Temlyakov.
We show the optimality of this bound with an explicit construction,
based on a refinement of the method of Olevskii. This construction leads to
other examples of quasi-greedy bases with large $k_N$ in Banach spaces, which are
of independent interest.
\end{abstract}

\section{Introduction }
\setcounter{equation}{0}\setcounter{footnote}{0}
\setcounter{figure}{0}

The concept of quasi-greedy basis evolved from the analysis of
thresholding algorithms  for non-linear $N$-term approximation in
Banach spaces; see e.g. \cite{Tem1} for a detailed presentation
and background. In recent years it has attracted attention from
both, the approximation theory and the Banach space point of view.

Let us recall the relevant definitions and standard notation. For
a (normalized) basis $\{\be_j\}_{j=1}^\infty$ in a Banach space
$\SX$ and $N=1,2,\dots$ we consider non-linear operators $G_N$ as
follows
\[ x=\sum_{j=1}^\infty a_j\be_j\in \SX\longmapsto
G_N(x)=\sum_{j\in \Lambda} a_j\be_j,\] where $\Lambda$ is {\em any
} $N$-element subset of $\{1,2,\dots\}$ such that $\min_{j\in
\Lambda}|a_j|\geq \max_{j\notin \Lambda}|a_j|$. Then $\{\be_j\}$
is called a {\em quasi-greedy basis} if for any $x\in \SX$ and any
choice of $G_N$'s we have $\lim_{N\to \infty}\|x-G_N(x)\|=0$,
 that is the series defining $x$ converges in norm after
decreasing rearrangement of their summands. It is known (see
\cite{Wo}) that this is equivalent to  \Be
\|G_Nx\|\leq K\|x\|, 
 \quad \forall\;x\in\SX,\;N=1,2,\ldots
\label{K}\Ee for some (smallest) constant $K$, which we assume
fixed throughout the paper. In particular, every unconditional
basis is quasi-greedy, but there exist also examples of
\emph{conditional} quasi-greedy bases \cite{KT,Wo, DKK,N1,Go,DST}.
In this paper we shall be interested in the latter.

\

Associated with a basis $\sB=\{\be_j\}$ in $\SX$, we consider the sequence
\[
k_N=k_N(\mathscr{B}):=\sup_{|A|\leq N}\,\big\|S_A\big\|,\quad
N=1,2,\ldots
\]
where $S_A:\SX\to\SX$ denotes the projection operator $S_A(x)=\sum_{j\in A}a_j(x)\be_j$.
 Generally speaking, the constants $k_N$ quantify
the conditionality of the basis $\sB$. In fact, $\sB$ is
unconditional if and only if $k_N(\sB)=O(1)$.

In approximation theory $k_N$ can also be used to quantify the
performance of greedy algorithms with respect to the best $N$-term
approximation from $\{\be_j\}$; that is, if $C_N$ denotes the
smallest constant such that
\[
\|x-G_Nx\|\,\leq \,C_N\,\inf\big\{\|x-\sum_{j\in A}c_j\be_j\|
\mid c_j\in\SC,\;|A|\leq N\,\big\},\quad\forall\;x\in\SX,
\]
then it is proved in \cite{GHO,TYY2} that $C_N\approx k_N$ when
 $\{\be_j\}$ is an \emph{almost-greedy} basis of $\SX$
(i.e. quasi-greedy and democratic\footnote{In Hilbert spaces,
quasi-greedy bases are always democratic \cite{Wo}, so both
concepts coincide.}).
Thus, in this case the
constants $k_N$ also give information on the rate of convergence
of greedy algorithms.

 \

It is known that for quasi-greedy bases in Banach spaces one has
\[k_N=O(\log N)\] (see \cite[Lemma 8.2]{DKK}), and this bound is
actually attained in some Banach spaces \cite{GHO}. It was asked
in \cite[p. 335]{TYY1} whether this bound is optimal or could be
improved in the case of Hilbert spaces. Our first result answers
this question.

\begin{theorem}
\label{Th1} Let $\SH$ be a Hilbert space and $\{\be_j\}$ a
quasi-greedy (normalized) basis with constant $K$. Then, there
exists $\al=\al(K)<1$ and $c>0$ such that \Be
k_N\big(\{\be_j\}\big)\leq\, c\, (\log N)^\al,\quad
\forall\;N\in\SN. \label{kNal}\Ee Moreover, if $\{\be_j\}$ is
besselian or hilbertian then one can choose $\al<\frac12$ in
\eqref{kNal}.
\end{theorem}

Recall that
 $\{\be_j\}$ is {\em besselian}
if $\sum_j|a_j|^2\leq C\,\|\sum_j a_j \be_j\|^2_\SX$ for all
finitely supported scalars
  $(a_j)$, and is called
{\em hilbertian }  if the converse inequality $ \sum_j|a_j|^2\geq
C\,\|\sum_j a_j \be_j\|^2_\SX$ holds.

\

Our second result proves that the bound obtained in \eqref{kNal}
is actually optimal.

\begin{theorem}
\label{Th2} For every $\al<1$, there exists a quasi-greedy basis
in
 $\SH$ and a constant $c_\al>0$ such that \Be k_N\geq c_\al(\log N)^{\al}, \quad N=1,2,\ldots\label{kNcal}\Ee
If $\al<1/2$, then the basis can be chosen to be in addition besselian (or hilbertian).
 \end{theorem}

 Theorem \ref{Th1} is shown in $\S2$, with
an explicit expression for $\al=\al(K)$ given in \eqref{alpha=}.
In the proof we make
use of the inner product structure of $\SH$, although the argument
can be adapted to other settings, such as $L^p$ spaces,
$1<p<\infty$, for which \eqref{kNal} is also true if $\{\be_j\}$
is quasi-greedy (see Appendix II).

Theorem \ref{Th2} is shown in $\S3$. The proof is based on a
construction due to Olevskii, which was developed in \cite{Wo} to
produce conditional quasi-greedy bases in Banach spaces. This
construction has an independent interest, and is stated as Theorem
\ref{ThO1} below. Its proof contains new ideas compared to
\cite[Theorem 2]{Wo}. Namely, we refine the method so that
besselian assumptions are not needed,
and moreover to obtain a basis which is almost-greedy and has
largest possible $k_N$. Only in this way we can reach the optimal
bounds in \eqref{kNcal}\footnote{The construction in \cite{Wo}
would only lead to $k_N\gtrsim(\log\log N)^{\al}$.}. We also apply
this construction to obtain new examples of almost-greedy bases in
Banach spaces with $k_N\approx\log N$.

\bline{\bf Acknowledgements:}  The authors thank Eugenio
Hern\'andez for useful conversations on this topic, and for
pointing out  a simplification in the original proof of Lemma
\ref{L1}. The second author also acknowledges the pleasant
atmosphere of the \emph{9th International Conference on Harmonic
Analysis (El Escorial 2012)}, where this research started.

\section{Proof of Theorem \ref{Th1}}

 Below we identify $x=\sum_{j=1}^\infty a_j(x)\be_j\in\SH$ with the coefficient sequence $(a_j)_{j=1}^\infty$,
so we write $\supp x=\{j\in\SN\mid a_j(x)\not=0\}$.
We shall use the following definition.

\begin{definition} Let $x,y\in\SH$. We say that $x\succcurlyeq y$ if
\Benu \item[(i)] $\supp x\cap\supp y=\emptyset$
\item[(ii)] $\min_{i\in\supp x}|a_i(x)|\geq\max_{j\in\supp y}|a_j(y)|$
\Eenu\end{definition}

The key result is the following lemma. Recall that the
quasi-greedy constant $K$ was defined in \eqref{K}.

\begin{lemma}\label{L1}
There exists $\dt=\dt(K)\in[0,1)$ such that, for all $x\succcurlyeq y$
\Be
(1-\dt)\Big(\|x\|^2+\|y\|^2\Big)\,\leq \big\|x+y\big\|^2\,\leq\,(1+\dt)\Big(\|x\|^2+\|y\|^2\Big).
\label{dt}\Ee
In fact, \eqref{dt} holds with $\dt=\sqrt{1-\frac1{K^2}}$.
\end{lemma}
\Proof
 Let $\ga\in\SR$ with $|\ga|\leq 1$. Then \eqref{K} implies
\[
\|x\|^2\leq K^2\big\|x-\ga y\big\|^2\,=\,K^2\big( \|x\|^2+\ga^2\|y\|^2-2\ga\Re\xy\big).\]
Hence we have the inequality
\[
\|y\|^2\,\ga^2-2\ga\Re\xy+\big(1-\tfrac1{K^2}\big)\|x\|^2\geq0,
\quad \forall\;\ga\in[-1,1].
\]
This inequality holds with $y$ replaced by $e^{i\theta}y$, for any
$\theta\in\SR$, so we also have \Be
\|y\|^2\,\ga^2-2\ga|\xy|+\big(1-\tfrac1{K^2}\big)\|x\|^2\geq0,
\quad \forall\;\ga\in[-1,1]. \label{aux3}\Ee Substituting
$\ga=\sqrt{1-\frac{1}{K^2}}$ into \eqref{aux3} we obtain \Be
2|\xy|\leq
\sqrt{1-\tfrac1{K^2}}\,\big(\|x\|^2+\|y\|^2),\label{2xy} \Ee from
which \eqref{dt} follows easily with $\dt=\sqrt{1-\tfrac1{K^2}}$.

\ProofEnd

As a special case of the lemma we obtain

\begin{corollary}
If $\{\be_j\}$ is a quasi-greedy (normalized) basis in $\SH$ such
that
\[
\|G_Nx\|\leq \|x\|,\quad\forall\;x\in\SH,\; N=1,2,\ldots
\]
then, $\{\be_j\}$ is an orthonormal basis.\end{corollary} \Proof
Applying Lemma \ref{L1}  to $x=\be_i$, $y=\be_j$ with $i\not=j$
and $K=1$, we obtain
 $\langle \be_i,\be_j\rangle=0$.
\ProofEnd

An iteration of the previous lemma leads to the following. We
denote, for $\ga\in\SR$, $\lceil \ga\rceil=\min\{k\in\SZ\mid
\ga\leq k\}$.

\begin{lemma}\label{L2}
Let $\dt=\dt(K)$ be as in Lemma \ref{L1}, then for all $x_1\succcurlyeq x_2\succcurlyeq \ldots \succcurlyeq x_m$
with pairwise disjoint supports
we have
\Be
(1-\dt)^{\lceil\log_2 m\rceil}\,\sum_{j=1}^m\|x_j\|^2\,\leq \big\|x_1+\ldots+x_m\big\|^2\,\leq\,(1+\dt)^{\lceil\log_2 m\rceil}\sum_{j=1}^m\|x_j\|^2.
\label{dtm}\Ee
\end{lemma}
\Proof We shall prove the result for $2^{n-1}<m\leq 2^n$ by
induction in $n=\lceil\log_2 m\rceil$. The case $n=1$ corresponds
to \eqref{dt}. Assume \eqref{dtm} holds for $m\leq 2^n$, and we
shall verify it for $2^{n}<m\leq 2^{n+1}$. Call $x'=\sum_{1\leq
j\leq 2^n} x_j$ and $y'=\sum_{2^n<j\leq m} x_j$. Since
$x'\succcurlyeq y'$, Lemma \ref{L1} gives\Beas
\big\|\sum_{j=1}^m x_j\big\|^2 =  \big\|x'+y'\big\|^2& \leq & (1+\dt)\Big(\|x'\|^2+\|y'\|^2\Big)\\
& = & (1+\dt)\Big(\|x_1+\ldots+x_{2^n}\|^2+\|x_{2^n+1}+\ldots+x_m\|^2\Big)\\
& \leq & (1+\dt)^{n+1}\sum_{j=1}^m\|x_j\|^2,
\Eeas
using  the induction hypothesis in the last step.
The inequality from below is similar.
\ProofEnd

From these two lemmas, the proof of Theorem \ref{Th1} is similar
to \cite[Theorem 5.1]{GHO} (see also \cite{DKK,DST}).

\

\Proofof{Theorem \ref{Th1}}

Let $A\subset\SN$ with $|A|=N\geq2$. We must show that, for all
$x=\sum_i a_i \be_i\in\SH$ then \Be \|S_A(x)\|\leq \,c\,(\log
N)^\al\,\|x\|, \label{aux0}\Ee for some $\al<1$ (independent of
$x$ and $N$).
 By scaling we may assume $\max_i |a_i| =
1$ (which by \eqref{K} implies
$ \|x\| \geq 1/K$).

Let $m=\lceil\log_2N\rceil$, so that $2^{m-1}<N\leq 2^m$.
For $\ell=1,\ldots,m$, we define
\[
F_\ell=\{j\mid 2^{-\ell}<|a_j|\leq 2^{-(\ell-1)}\}
\mand F_{m+1}=\{j\mid |a_j|\leq 2^{-m}\}.\]
Next write $A$ as a disjoint union of the sets
$A_\ell=A\cap F_\ell$, {\small $\ell=1,\ldots, m+1$}.
Clearly
 \Be\|S_{A_{m+1}}x\| \leq \sum_{i\in A_{m+1}} |a_i| \|\be_i\| \leq
 2^{-m}N \leq 1 \leq K
 \|x\|.\label{Am1}\Ee
For the other terms we shall use Lemmas 5.2 and 5.3 in \cite{GHO}, which give
 \[ \|S_{A_\ell} x\| \leq \,c_1\, \|S_{F_\ell}(x)\| \leq \,c_2\, \|x\|,\]
with $c_1=64K^3$  and $c_2=128K^4$.  Now, Lemma \ref{L2} gives
 \Be\|\sum_{\ell=1}^m S_{A_\ell} x\|^2  \leq    (1+\dt)^{\lceil\log_2 m\rceil}\sum_{\ell=1}^m \|S_{A_\ell} x\|^2
  \leq c_1^2(1+\dt)^{\lceil\log_2 m\rceil}\sum_{\ell=1}^m \|S_{F_\ell} x\|^2.\label{aux1}\Ee
We now have two possible approaches. In  the first approach we use the lower bound in \eqref{dtm},
so that \eqref{aux1} becomes\[
\|\sum_{\ell=1}^m S_{A_\ell} x\|^2  \leq   c_1^2 \Big(\tfrac{1+\dt}{1-\dt}\Big)^{\lceil\log_2 m\rceil}
\big\|\sum_{\ell=1}^m S_{F_\ell} x\big\|^2\,\leq\, c_3^2\Big(\tfrac{1+\dt}{1-\dt}\Big)^{\lceil\log_2 m\rceil}\,\|x\|^2,
\]
with $c_3=Kc_1=64K^4$. Observe that
\[
\Big(\tfrac{1+\dt}{1-\dt}\Big)^{\log_2 m}= 2^{\log_2 m\,\log_2\tfrac{1+\dt}{1-\dt}}=m^{\log_2\tfrac{1+\dt}{1-\dt}}=m^{2\al_1}
\]
if we take $\al_1=\frac12\log_2\tfrac{1+\dt}{1-\dt}$. Notice however that $\al_1<1$ if and only if $\dt<3/5$, so
this approach is not good for $\dt$ close to 1 (ie, when $K$ is very large).

\

A second approach for \eqref{aux1} consists in estimating each
$\|S_{A_\ell} x\|\leq c_2\|x\|$. Then \Be \|\sum_{\ell=1}^m
S_{A_\ell} x\|^2  \leq   c_2^2 (1+\dt)^{\lceil\log_2 m\rceil}
\,m\,\|x\|^2. \label{aux1bis}\Ee Now we can write\[(1+\dt)^{\log_2
m} \,m= 2^{\log_2
m\,\log_2(1+\dt)}\,m=m^{1+\log_2(1+\dt)}=m^{2\al_2},
\]
if we choose $\al_2=(1+\log_2(1+\dt))/2$. Notice that this time
$\al_2<1$, but it may happen that $\al_2>\al_1$ if $\dt< 1/2$. In
that case (ie, when $K$ is close to 1) the former choice is
slightly better.

Combining the two approaches, and using also \eqref{Am1}, we see
that \eqref{kNal} holds with
\Be
\al=\frac12\,\min\Big\{\log_2\tfrac{1+\dt}{1-\dt},\, 1+\log_2(1+\dt)\,\Big\},
\label{alpha=}\Ee
with the minimum attained in the first number for $\dt\leq\frac12$, and in the second number for $\dt\geq\frac12$.
Recall also from Lemma \ref{L1} that $\dt=\sqrt{1-1/K^2}$.

\

Finally, suppose that the basis $\{\be_j\}$ is not only quasi-greedy, but also
\emph{besselian}.
Quasi-greediness implies that $\ell^{2,1}\hookrightarrow\SH$ (\cite[Thm 3]{Wo}), so we can estimate for each $\ell$
\[
\|S_{F_\ell} x\| \lesssim 2^{-\ell}\,|F_\ell|^{1/2}\,\leq \big(\sum_{j\in F_\ell}|a_j|^2\big)^{1/2}.
\]
Inserting this into \eqref{aux1} and using that $\SH\hookrightarrow\ell^2$ (from the besselian assumption), we obtain
\[
\|\sum_{\ell=1}^m S_{A_\ell} x\|^2  \lesssim (1+\dt)^{\lceil\log_2 m\rceil}\,\sum_j|a_j|^2\lesssim m^{\log_2(1+\dt)}\,\|x\|^2.\]
Thus, \eqref{kNal} holds with $\al=\frac12\log_2(1+\dt)$, which is always a real number $<1/2$.
The same bound holds when $\{\be_j\}$ is \emph{hilbertian}, since in this case the dual basis $\{\be^*_j\}$
is besselian in $\SH^*$ (and also quasi-greedy, by \cite{DKKT}), while $k_N$ is the same for both bases.
\ProofEnd

\

\section{The proof of Theorem \ref{Th2}}
\subsection{A general construction of quasi-greedy bases}

 The next result gives a general method to produce
quasi-greedy bases in Banach spaces. As mentioned in the
introduction, it is also an improvement over the statement in
\cite[Theorem 2]{Wo}.

\begin{theorem}
 \label{ThO1}
Let $\cX$ be a Banach space with a (normalized) basis
$\mathfrak{X}=\{x_k\}_{k=1}^\infty$. Then, the
space\footnote{Endowed with the norm
$\|x+y\|_{\cX\oplus\ell^2}^2=\|x\|^2_\cX+\|y\|^2_{\ell^2}$.}
$\cX\oplus \ell^2$ has a quasi-greedy basis $\Psi$.

 Moreover, \Benu \item[(i)] $\Psi$ is democratic and
 $\|\sum_{\la\in\La}\psi_\la\|\approx|\La|^{1/2}$.

 \item [(ii)] if the basis $\fX$ is
besselian (or hilbertian), so is $\Psi$. \item[(iii)] if the basis
$\fX$ has the property that, for some $c>0$ and every
$N=1,2,\ldots$
 \Be \exists\;x\in\cX \mbox{ and
}A\subset\{1,\ldots,N\} \quad\mbox{such that}\quad\|S_{A}x\|\geq
\,c\,{k_N(\fX)}\,\|x\|, \label{xA}\Ee then the quasi-greedy basis
$\Psi$ satisfies
\[ k_N(\Psi)\gtrsim k_{\log N}(\fX),\quad N=2,3,\ldots\]
\Eenu
\end{theorem}

 \Proof Write $\mathfrak{X}=\{x_k\}_{k=1}^\infty$ for
the basis in $\cX$, and $\{e_j\}_{j=1}^\infty$ for the canonical
orthonormal basis in $\ell^2$. In the direct sum
 $\cX\oplus\ell^2$, consider the system of vectors $\Upsilon$ given by\Be x_1,e_1;\;
x_2,e_2,\ldots, e_{n_3};\,x_3, e_{n_3+1},\ldots, e_{n_4}; \ldots
\label{xe}\Ee for a suitable increasing sequence $n_k$. Here we
choose $n_1=0$, $n_2=1$ and $n_{k+1}=n_k+2^k-1$, so that each
block $\Upsilon_k=\{x_k,e_{n_k+1},\ldots, e_{n_{k+1}}\}$ generates
a subspace $\SH_k$ of dimension $N_k:=n_{k+1}-n_k+1=2^k$, of which
$\Upsilon_k$ is a natural orthonormal basis. We rename this basis
as
\[ \Upsilon_k=\big\{g_{k,1},\ldots, g_{k,2^k}\big\},\] and write
the system in \eqref{xe} as $\Upsilon=\cup_{k=1}^\infty
\Upsilon_k$. The next lemma follows from elementary Banach space
theory.

\begin{lemma}\label{3.1}
The system $\Upsilon$ in \eqref{xe} is a basis in
$\cX\oplus\ell^2$. Moreover, if $\mathfrak{X}$ is besselian (or
hilbertian), so is $\Upsilon$.\end{lemma}

We now use the Olevskii construction; see \cite{Wo}. For each $k$, let $A=A^{(k)}$ denote the
matrix in $SO(2^k,\SR)$ with entries given by the Haar basis in $\SR^{2^k}$, ie
\[
A^{(k)}=\rm{col}\;\big\{h_0,h_1,\ldots, h_{2^k-1}\big\}
\]
 where $h_0=2^{-k/2}$. In each subspace $\SH_k$ we define a new orthonormal basis
$\Psi_k=\{\psi_{k,1},\ldots,\psi_{k,2^k}\}$,  by letting \Be
 \left(
   \begin{array}{c}
     \psi_{k,1}\\
     \vdots \\
     \psi_{k,2^k} \\
   \end{array}
 \right)\,=\, \Bigg( \;A^{(k)}\;\Bigg) \left(
   \begin{array}{c}
     g_{k,1}\\
     \vdots \\
     g_{k,2^k} \\
   \end{array}
 \right)\label{matrix}\Ee
 That is,
 \Be
 \psi_{k,\ell}\,=\,2^{-k/2} x_k\,+\,\sum_{m=2}^{2^k} a_{\ell,m}\,g_{k,m},\quad \ell=1,\ldots, 2^k.
 \label{aux3.1}\Ee
 From the orthonormality of $\Psi_k$ and Lemma \ref{3.1} it easily follows that

 \begin{lemma}\label{3.2}
 The system $\Psi=\cup_{k=1}^\infty\Psi_k$ is a basis of $\cX\oplus\ell^2$.
 Again, if $\mathfrak{X}$ is besselian (or hilbertian), so is $\Psi$.\end{lemma}

The key step of Theorem \ref{ThO1} is to establish the
quasi-greediness of $\Psi$.
 For this we need to refine the analysis of  Olevskii construction given in \cite{Wo}.
 Notice that we do not require the basis $\mathfrak{X}$  to be besselian in $\cX$.

\begin{lemma}\label{3.3}
  $\Psi$ is a quasi-greedy  basis of $\;\cX\oplus\ell^2$, that is
  \Be
  \|G_N(z)\|\leq\,C\,\|z\|,\quad \forall\;z\in\cX\oplus\ell^2,\quad N=1,2,\ldots
  \label{Gnz}\Ee\end{lemma}
  \Proof
  Let $z=\sum_{k=1}^\infty\sum_{\ell=1}^{2^k} c_{k,\ell}\psi_{k,\ell}$, and  $\La:=\supp G_N(z)$.
  We use the notation $P_\cX, P_\cH$ for the natural projections onto $\cX$ and $\cH=\ell^2$
  respectively, and $S_\La$ for the projection onto $\span\{\psi_\la\}_{\la\in\La}$.
  We also write $\La_k=\{\ell\mid (k,\ell)\in\La\}$.

  We need to show that\[
  \|S_\La z\|^2=\,\|P_\cX S_\La  z\|^2+\|P_\cH S_\La  z\|^2\,\leq C\,\|z\|^2.\]
  We begin with the first summand on the left hand side; that is, we shall show that
\Be \|P_\cX S_\La  z\|\,\leq C\,\|z\|.\label{aux3.2}\Ee Let
$\al=\min_{(k,\ell)\in\La}|c_{k,\ell}|$, which we may assume
$\al>0$ (otherwise $G_Nz=z$ and \eqref{Gnz} is trivial). Fix
$M\geq1$ to be chosen later, and notice from \eqref{aux3.1} that
we can split \Be P_\cX S_\La (z) \,=\,\sum_{k\geq
M}2^{-k/2}\big(\sum_{\ell\in\La_k} c_{k,\ell}\big)
x_k\,+\,\sum_{k< M}2^{-k/2}\big(\sum_{\ell\in\La_k}
c_{k,\ell}\big) x_k. \label{PXsplit}\Ee The first term has norm
bounded by \Bea \big\| \sum_{k\geq
M}2^{-k/2}\big(\sum_{\ell\in\La_k} c_{k,\ell}\big) x_k\big\| &
\leq & \sum_{k\geq
M}2^{-k/2}\sum_{\ell\in\La_k}|c_{k,\ell}|\nonumber\\
& \leq & \sum_{k\geq M}2^{-k/2}
\,\frac1\al\,\sum_{\ell\in\La_k}|c_{k,\ell}|^2\nonumber\\
& \leq& c\,2^{-M/2}\,\frac1{\al}\,\sup_k\|P_{\SH_k}z\|^2\,\leq\,\,c'\,2^{-M/2}\,\frac{\|z\|^2}{\al}.\label{aux3.4}\Eea
If $\|z\|\leq 2\al$, we can choose $M=1$ and we are done. Otherwise
\[
\big\| \sum_{k< M}2^{-k/2}\big(\sum_{\ell\in\La_k} c_{k,\ell}\big) x_k\big\| \leq
\big\| \sum_{k< M}2^{-k/2}\big(\sum_{\ell=1}^{2^k} c_{k,\ell}\big) x_k\big\|  +
\big\| \sum_{k< M}2^{-k/2}\big(\sum_{\ell\not\in\La_k} c_{k,\ell}\big) x_k\big\| \,=\,I\,+\,II.
\]
Clearly,
\[
I=\big\|P_\cX\big( S_{\{k<M\}}z\big)\big\|\leq
\big\|S_{\{k<M\}}(z)\big\|\leq\,C\,\|z\|\] since $\Psi$ is a
basis. Finally, since $\max_{(k,\ell)\not\in\La}|c_{k,\ell}|\leq
\al$, \Be II\leq \sum_{k<M}
2^{-k/2}\sum_{\ell\not\in\La_k}|c_{k,\ell}|\leq \al\sum_{k<M}
2^{k/2}\leq c\, \al\, 2^{M/2}.\label{aux3.5}\Ee So we can optimize
in \eqref{aux3.4} and \eqref{aux3.5} by choosing $M$ such that
$2^M=\|z\|^2/\al^2$. This proves \eqref{aux3.2}.

\medskip

Next we show that \Be \|P_\cH S_\La  z\|\,\leq
C\,\|z\|.\label{auxPHz}\Ee This would be easy to establish if we
assume that $\mathfrak{X}$ is besselian. Indeed, in that case,
using the orthogonality of the spaces $P_\cH(\SH_k)$ we can write
\Bea
\big\|P_\cH S_\La z\big\|^2 
 & = & \sum_k \big\|\sum_{\ell\in\La_k} c_{k,\ell} P_\cH(\psi_{k,\ell})\big\|^2\nonumber\\
   & \leq & \sum_k \big\|\sum_{\ell\in\La_k} c_{k,\ell} \psi_{k,\ell}\big\|^2\,=\,\sum_k\sum_{\ell\in\La_k}|c_{k,\ell}|^2\nonumber\\
  & \leq &  \sum_k\sum_{\ell=1}^{2^k}|c_{k,\ell}|^2\,\leq\, C\,\|z\|^2,\label{PHs}\Eea
where in the last inequality we would use that the basis $\Psi$ is
also besselian.

\

We now give a different argument which holds for general $\fX$.
As before, we define $\al=\min_{(k,\ell)\in\La}|c_{k,\ell}|$ which we may assume $\al>0$.
We write $z_k=P_{\SH_k}(z)=\sum_{\ell=1}^{2^k}c_{k,\ell}\psi_{k,\ell}$ as
\[z_k=P_\cX(z_k)+P_\cH(z_k)=\la_k x_k+\sum_{\ell=1}^{2^k}\eta_{k,\ell}\psi_{k,\ell},\]
for suitable scalars $\la_k$ and $\eta_{k,\ell}$. Since from \eqref{matrix} we have
$x_k=2^{-k/2}\sum_{\ell=1}^{2^k}\psi_{k,\ell}$, we see that
\Be
c_{k,\ell}=2^{-k/2}\la_k\,+\,\eta_{k,\ell},\quad \ell=1,\ldots,2^k.
\label{ckl}\Ee
We want to show that
\Be
\big\|P_\cH(S_\La z)\big\|^2\,=\,\sum_{k=1}^\infty\big\|P_\cH(S_\La z_k)\big\|^2\,\lesssim\,
\|z\|^2+\,\sum_{k=1}^\infty\big\|P_\cH(z_k)\big\|^2,\label{C0}
\Ee
from which \eqref{auxPHz} would follow easily, since the last series equals $\|P_\cH(z)\|^2\leq\|z\|^2$.
To establish \eqref{C0} we consider three possible situations for the index $k$,
\begin{eqnarray*}
 A_1&=&\big\{ k \ :\ 2^{-k/2}|\lambda_k|<\al/2\big\}\\
A_2&=&\big\{ k\ :\ 2^{-k/2}|\lambda_k|\in[\tfrac\al2,2\al]\big\}\\
A_3&=&\big\{ k \ :\ 2^{-k/2}|\lambda_k|>2\al \big\}.
\end{eqnarray*}
Assume that $k\in A_1$. Then
\[
\big\|P_\cH (S_\La z_k)\big\|^2 \,\leq\, \big\|\sum_{\ell\in\La_k}
c_{k,\ell}
\psi_{k,\ell}\big\|^2\,=\,\sum_{\ell\in\La_k}|c_{k,\ell}|^2.
\]
Now, when $\ell\in \La_k$ we have $|c_{k,\ell}|\geq \al$, and hence by \eqref{ckl}
\[
\al\leq \big|2^{-k/2}\la_k+\eta_{k,\ell}\big|\leq\,\frac\al2+|\eta_{k,\ell}|.
\]
Thus, $\tfrac\al2\leq|\eta_{k,\ell}|$, which in turn implies
$|c_{k,\ell}|\leq 2|\eta_{k,\ell}|$. We conclude that, for $k\in
A_1$, \Be \big\|P_\cH(S_\La
z_k)\big\|^2\,\leq\,\sum_{\ell\in\La_k}|c_{k,\ell}|^2\,\leq\,4\,\sum_{\ell=1}^{2^k}|\eta_{k,\ell}|^2\,=\,4\,\big\|P_\cH(z_k)\big\|^2.
\label{C1}\Ee Next we consider $k\in A_2$. Here we use the cruder
bound $\|P_\cH(S_\La z_k)\|\leq \|S_\La z_k\|\leq\|z_k\|$, and
notice that
\[
\sum_{k\in A_2}\|z_k\|^2\,=\,\sum_{k\in
A_2}\Big(\|P_\cX(z_k)\|^2+\|P_\cH(z_k)\|^2\Big)= \sum_{k\in
A_2}|\la_k|^2+\sum_{k\in A_2}\|P_\cH(z_k)\|^2.
\]
We thus need to bound $\sum_{k\in A_2}|\la_k|^2$. Notice that $A_2$ is a finite set (since $2^{-k/2}\la_k\to0$ as $k\to\infty$), and write
$N_0=\max A_2$. Clearly,
\[
2^{-N_0/2}|\la_{N_0}|\,\approx\,\al.
\]
Then,
\[
\sum_{k\in A_2}|\la_k|^2\leq 4\al^2\sum_{k\leq N_0} 2^k\leq 8\al^2 2^{N_0}\lesssim \,|\la_{N_0}|^2.
\]
Since $|\la_{N_0}|=\|P_\cX(z_{N_0})\|\leq C\|z\|$, we see that \Be
\sum_{k\in A_2} \|P_\cH(S_\La
z_k)\|^2\,\lesssim\,\|z\|^2\,+\,\sum_{k\in A_2}\|P_\cH(z_k)\|^2.
\label{C2}\Ee

Finally, consider $k\in A_3$. Using once again \eqref{ckl} we see that\[
S_\La z_k=2^{-\frac k2}\la_k\sum_{\ell\in\La_k}\psi_{k,\ell}\,+\,\sum_{\ell\in\La_k}\eta_{k,\ell}\psi_{k,\ell},
\]
and therefore\Be
\big\|P_\cH(S_\La z_k)\big\|\leq\,2^{-\frac k2}|\la_k|\,\big\|P_\cH(\sum_{\ell\in\La_k}\psi_{k,\ell})\big\|\,+\,
\Big(\sum_{\ell=1}^{2^k}|\eta_{k,\ell}|^2\Big)^\frac12.
\label{auxA3}\Ee
The second summand equals $\|P_\cH(z_k)\|$, so we will work on the first.

 Notice that $P_\cH(\sum_{\ell=1}^{2^k}\psi_{k,\ell})=P_\cH(2^{k/2}x_k)=0$, so we
 have\Be
\big\|P_\cH(\sum_{\ell\in\La_k}\psi_{k,\ell})\big\|^2=
\big\|P_\cH(\sum_{\ell\in\La_k^c}\psi_{k,\ell})\big\|^2\leq|\La_k^c|.
\label{Gamma}\Ee Now, $k\in A_3$ and $\ell\in \La_k^c$ imply
$|c_{k,\ell}|\leq \al$. Using \eqref{ckl} we see that
\[
2^{-\frac k2}|\la_k|-|\eta_{k,\ell}|\,\leq\, \big|2^{-\frac k2}\la_k+\eta_{k,\ell}| \,\leq \,\al\,\leq \frac{2^{-\frac k2}|\la_k|}2,
\]
and hence $|\eta_{k,\ell}|\geq2^{-k/2}|\la_k|/2$. Therefore, using
also \eqref{Gamma}, the middle term in \eqref{auxA3} is bounded by
\[
2^{-k}|\la_k|^2|\La_k^c|\,\leq\,4\,\sum_{\ell\in\La_k^c}|\eta_{k,\ell}|^2\,\leq\,4\,\sum_{\ell=1}^{2^k}|\eta_{k,\ell}|^2
=4\big\|P_\cH(z_k)\|^2.\] Thus, for $k\in A_3$ in also have \Be
\big\|P_\cH(S_\La z_k)\big\|\leq\,3\,\big\|P_\cH(z_k)\big\|.
\label{C3} \Ee Thus, we can now combine \eqref{C1}, \eqref{C2} and
\eqref{C3} to obtain the asserted estimate in \eqref{C0}, and
hence establish Lemma \ref{3.3}. \ProofEnd

\begin{lemma}
\label{demo} The basis $\Psi$ is democratic and, for every finite
$\La$,
\Be\|\sum_{\la\in\La}\psi_{\la}\|\approx|\La|^{1/2}.\label{demo2}\Ee
\end{lemma}
\Proof The proof is a small refinement of the previous arguments.
For simplicity, we use the notation
$\bone_\La=\sum_{\la\in\La}\psi_{\la}$, and set
$\La_k=\{\la\in\La\mid \la=(k,\ell)\mbox{ for some $\ell$}\}$.
Call $N=|\La|$ and $N_k=|\La_k|$.

We first find an upper bound for \Be
\|\bone_\La\|^2=\|P_\cX(\bone_\La)\|^2+\|P_\cH(\bone_\La)\|^2.\label{bone}\Ee
Arguing as in \eqref{PHs}, the second term is easily estimated
by\[ \big\|P_\cH \bone_\La \big\|^2  = \sum_k
\big\|P_\cH(\bone_{\La_k})\big\|^2
 \leq  \sum_k \big\|\bone_{\La_k}\big\|^2\,=\,\sum_kN_k=N.\]
For the first, since $N_k\leq \min\{2^k,N\}$, setting $M=\log_2N$,
and arguing as in \eqref{PXsplit} \[ \big\|P_\cX \bone_\La \big\|
= \big\|\sum_{k=1}^\infty 2^{-k/2}N_kx_k\big\| \, \leq \,
\sum_{1\leq k\leq M}
2^{k/2}\,+\,\sum_{k>M}2^{-k/2}N\,\leq\,c\,\sqrt N.
\]
We now find a lower bound for \eqref{bone}. Partition the indices
$k$ by
\[
I_0=\{k\mid N_k\leq 2^{k-1}\}\mand I_1=\{k\mid 2^{k-1}<N_k\leq
2^k\}.
\]
Let $k_1=\max I_1$. Then, since $N_k\approx 2^k$ for $k\in I_1$,
we have \Be \big\|P_\cX\bone_\La\big\|^2=\big\|\sum_{k=1}^\infty
2^{-k/2}N_kx_k\big\|^2\geq c\,2^{-k_1}N^2_{k_1}\geq c' \sum_{k\in
I_1} N_k. \label{PXk}\Ee For the other term we use the identity
\Be
\big\|P_\cH\bone_{\La_k}\big\|^2=N_k(1-2^{-k}N_k).\label{bonek}
\Ee Assuming \eqref{bonek}, one sees that\[ \big\|P_\cH \bone_\La
\big\|^2  \geq \sum_{k\in I_0}
\big\|P_\cH(\bone_{\La_k})\big\|^2\geq \frac12\sum_{k\in I_0}N_k,
\]
which combined with \eqref{PXk} gives $\|\bone_\La\|^2\gtrsim N$.
It remains to show \eqref{bonek}, but this is easy, since by
orthogonality\[
N_k=\|\bone_{\La_k}\|^2=\|P_\cX\bone_{\La_k}\|^2+\|P_\cH\bone_{\La_k}\|^2
= 2^{-k}N_k^2+\|P_\cH\bone_{\La_k}\|^2,
\]
from which the claim follows easily.
 \ProofEnd

\

Finally we give a bound for $k_N(\Psi)$ in terms of $k_N(\fX)$.
\begin{lemma}\label{knH}
Assume that the basis $\fX=\{x_n\}$ of $\cX$ satisfies the
property in \eqref{xA}. Then, the basis $\Psi$ of
$\cX\oplus\ell^2$ constructed above has\Be k_M(\Psi)\,\geq \,c'\,
k_{\log M}(\fX),\quad M=2,3,\ldots \label{knKH}\Ee
\end{lemma}
\Proof Fix $M$ and choose $N$ such that $2^{N-1}\leq M<2^{N}$.
Select $x$ and $A$ as in \eqref{xA}, and set $\La= \cup_{k\in
A}\{\psi_{k,1},\ldots,\psi_{k,2^k}\}$, which has cardinality
$|\La|\leq 2^{N+1}$. Then
\[
\big\|S_\La x\big\|\geq\big\|P_\cX S_\La x\big\|=\big\|S_A
x\big\|\geq \,c\,{k_N(\fX)}\,\|x\|\geq \,c \,k_{\log_2
M}(\fX)\,\|x\|.\] Since $k_N$ is doubling, \eqref{knKH} follows
easily. \ProofEnd

\medskip
This completes the proof of Theorem \ref{ThO1}.\ProofEnd

\subsection{Conditional bases with large $k_N$}

We shall use Theorem \ref{ThO1} to prove Theorem \ref{Th2}. For
this purpose we first need to find a Hilbert space $\cX$ with a
conditional basis $\fX$ (not necessarily quasi-greedy) having
$k_N(\fX)$ as large as possible.
Here we give two examples in this direction.

We need the following elementary lemma about the Dirichlet
kernel $D_N(t)=\sum_{|n|\leq N}
e^{int}=\sin(N+\frac12)t/\sin(\frac t2)$.

\begin{lemma}
\label{dir} Let $|\ga|<1$. Then, \Be \label{dir1}
\big\|D_N(t)\big\|_{L^2(|t|^\ga\,dt)} \,\approx\,
N^\frac{1-\ga}2,\quad N=1,2,\ldots\Ee
\end{lemma}
\Proof The proof is elementary. From below,
\[
\big\|D_N(t)\big\|_{L^2(|t|^\ga\,dt)}^2\,\geq\, \int_{|t|\leq
1/N}\,|D_N(t)|^2\,|t|^{\ga}\,dt\,\approx\,
 N^2\int_{|t|\leq 1/N}\,|t|^{\ga}\,dt\approx
 N^{1-\ga}.
\]
From above, the remaining part of the integral is estimated by
\[ \int_{\frac1N\leq|t|<\pi}|D_N(t)|^2|t|^\ga\,dt\lesssim
\int_{\frac1N}^\pi t^{\ga-2}\,dt\,\lesssim \,N^{1-\ga}.
\] \ProofEnd

In the first example we obtain a \emph{besselian} basis with
$k_N\gtrsim N^{\frac12-\e}$.

\begin{proposition}
\label{P1} If $\al\in(0,1/2)$ then  there is a Hilbert space $\cX$
with a besselian conditional basis $\fX=\{x_n\}_{n=1}^\infty$ such
that \[k_N(\fX)\gtrsim N^\al.\] Moreover, for every $N=1,2,\ldots$
 there exists a partition $\{1,\ldots,N\}=A\uplus B$ such that
 \Be
 \big\|\sum_{n\in A}x_n\big\|_\cX\,\gtrsim\, N^{\al}\,\big\|\sum_{n\in A}x_n-\sum_{n\in B}x_n\big\|_\cX.
 \label{N12a}\Ee
 \end{proposition}
 \Proof
 We consider the example proposed by Babenko \cite{Ba}. That is, we set  $\cX= L^2([-\pi,\pi],|t|^{-2\al}dt)$
 with the  usual trigonometric system $\fX=\{1,e^{it},e^{-it},e^{2it},e^{-2it},\ldots\}$.
 That $\fX$ is a basis follows from the fact that $|t|^{-2\al}$ is an $A_2$-weight when $|\al|<1/2$
 (see e.g. \cite{duo}).
 When $\al>0$ the weight is bounded from below by a positive constant, so
 $\|\sum a_ne^{int}\|_\cX\geq c\,\|\sum a_ne^{int}\|_{L^2}=c(\sum|a_n|^2)^{1/2}$,
 and the basis is besselian.

 We now prove \eqref{N12a}. By the lemma
 \[
 \big\|\sum_{|n|\leq N} e^{int}\big\|_\cX^2\, \gtrsim\,
 N^{1+2\al}.\]
 On the other hand, Khintchine's inequality gives
 \[
 \mathbb{E}\Big[\int_{\ST}\,\big|\sum_{|n|\leq N}\pm e^{int}\big|^2\,|t|^{-2\al}\,dt\Big]
\,\approx\,\int_{\ST}\,N\,|t|^{-2\al}\,dt\,\approx\, N,
 \]
 so for a certain fixed constant $C>0$, there must exist some choice of signs $\pm1$ such that $\|\sum_{|n|\leq N}\pm e^{int}\|_\cX\leq C\,
N^{1/2}$. Partitioning  $\{-N,\ldots,N\}=A_+\uplus A_-$ according to these signs,
 we have shown that
 \[
 \big\|\sum_{|n|\leq N} e^{int}\big\|_\cX \,\gtrsim\, N^\al\,\|\sum_{n\in A_+} e^{int}-\sum_{n\in A_-} e^{int}\|_\cX.
 \]
 Using that $\|\sum_{|n|\leq N} e^{int}\|_\cX\leq 2\|\sum_{n\in A} e^{int}\|_\cX$ either for $A=A_+$ or $A=A_-$,   \eqref{N12a} follows
easily.

 \ProofEnd

\BR\label{R1} The example in Proposition \ref{P1} actually
satisfies $k_N\approx N^\al$.  We sketch a proof of the upper
bound in Appendix I. \ER

In our second example we find a basis in a Hilbert space with
$k_N\gtrsim N^{1-\e}$. It is a consequence of a well-known theorem of Gurarii and
Gurarii \cite{GG} that this growth is best possible.

\begin{proposition}
\label{P1bis} Let $\al<1$. Then there is a Hilbert space $\cX$
with a conditional basis $\fX=\{x_n\}_{n=1}^\infty$ such that
$k_N(\fX)\gtrsim N^\al$. Moreover, there exists $c>0$ such that
for every $N=1,2,\ldots$
 there is a set  $A\subset\{1,\ldots,N\}$ and a non-null $x\in\cX$  so that
 \Be
 \big\|S_A(x)\big\|_\cX\,\geq\,c\, N^{\al}\,\big\|x\big\|_\cX.
 \label{N1a}\Ee
 \end{proposition}
\Proof Consider $\cX=L^2(\ST, |t|^\al dt)\oplus L^2(\ST,
|t|^{-\al} dt)$ for the $\ell^2$-direct sum. Call
$\{e_k\}_{k=1}^\infty$ and $\{f_k\}_{k=1}^\infty$ the respective
trigonometric bases as in the previous proposition, and consider a
new basis in $\cX$ given by\Be
x_{2k-1}=\frac{e_k+f_k}{\sqrt2},\quad
x_{2k}=\frac{e_k-f_k}{\sqrt2},\quad k=1,2,\ldots \label{xef}\Ee
For $N\geq1$, let $x=\sum_{n=1}^{2N}x_n=\sqrt2\sum_{k=1}^N e_k$.
If $N=2m+1$, then $x$ coincides with the Dirichlet kernel
$D_m(t)$, which by Lemma \ref{dir} gives the estimate
\[
\|x\|_\cX\,=\,\sqrt2\,\big\|D_m(t)\big\|_{L^2(|t|^\al\,dt)}
\,
\approx N^\frac{1-\al}2.
\]
Now let $A=\{1,3,\ldots,2N-1\}$, so that
$S_A(x)=\frac1{\sqrt2}\sum_{k=1}^Ne_k+\frac1{\sqrt2}\sum_{k=1}^Nf_k$,
which has norm
\[
\big\|S_A(x)\big\|_\cX\,\geq\,\tfrac1{\sqrt2}\,\big\|\sum_{k=1}^N
f_k\big\|\,=\,\tfrac1{\sqrt2}\,\big\|D_m(t)\big\|_{L^2(|t|^{-\al}\,dt)}
\,\approx N^\frac{1+\al}2.
\]
Combining these two estimates we obtain
$k_{N}\,\geq\,{\|S_A(x)\|}/{\|x\|}\,\gtrsim\, N^\al$, as well as
the assertion in \eqref{N1a}.
 \ProofEnd

\subsection{End of the proof of Theorem \ref{Th2}} Fix $\al<1$ and apply Theorem
\ref{ThO1} to the Hilbert space $\cX$ and the basis $\fX=\{x_n\}$
in Proposition \ref{P1bis}. This produces a quasi-greedy basis
$\Psi$ in the Hilbert space $\cX\oplus\ell^2$. Since  the
assumptions of Lemma \ref{knH} hold with $k_N(\fX)\approx N^\al$
(by \eqref{N1a}), we obtain that
\[
k_M(\Psi)\geq c_\al\,(\log M)^\al,\quad M=1,2,\ldots
\]
as we wished to prove.

If we only assume $\al<1/2$, then we would argue similarly, using
instead Proposition \ref{P1}, so that the basis $\Psi$ of
$\cX\oplus\ell^2$ is in addition besselian. In this case, the dual
system to $\Psi$ will be a hilbertian quasi-greedy basis for
$(\cX\oplus\ell^2)^*$ with the same bound on $k_M$. \ProofEnd \


\subsection{Further results}

 As a consequence of Theorem
\ref{ThO1} we can find new examples of quasi-greedy democratic
bases in Banach spaces for which $k_N\approx\log N$ (see also \cite{GHO}).

\begin{corollary}
\label{col5.2}There exists a basis in $c_0\oplus
\ell^1\oplus\ell^2$ which is quasi-greedy, democratic and has
$k_N\approx \log N$. \end{corollary} \Proof Consider $c_0$ and
$\ell^1$ with their usual canonical bases, say $\{e_k\}$ and
$\{f_k\}$, and define in $\cX=c_0\oplus \ell^1$ a basis
$\fX=\{x_n\}$ as in \eqref{xef}. Then a similar proof as in
Proposition \ref{P1bis} shows that $k_N(\fX)\gtrsim N$ and
\eqref{N1a} holds with $\al=1$. Then, Theorem \ref{ThO1} gives a
quasi-greedy democratic basis $\Psi$ in
$\cX\oplus\ell^2=c_0\oplus\ell^1\oplus\ell^2$ with
$k_N(\Psi)\gtrsim\log N$.\ProofEnd

Theorem \ref{ThO1} can also be used to show the optimality of the bound $k_N=O(\log N)^\al$
with $\al<1$ for quasi-greedy bases in $\ell^p$ or $L^p$ spaces, $1<p<\infty$ (see Theorem \ref{Thp} below).

\begin{corollary}
\label{collp}Let $1<p<\infty$ and $\al<1$. Then, there exists a
quasi-greedy basis in $\ell^p$ with $k_N\gtrsim(\log N)^\al$.
\end{corollary}
\Proof Use a similar construction to \cite[Corollary 5]{Wo}.
\ProofEnd

In the spaces $L^p[0,1]$ more can be said. Namely, a direct application of Theorem \ref{ThO1} gives bases which are additionally democratic,
with democracy function $N^{1/2}$.

\begin{corollary}
\label{colLp}Let $1<p<\infty$ and $\al<1$. Then, there exists a
quasi-greedy democratic basis $\Psi$ in $L^p[0,1]$ with $k_N(\Psi)\gtrsim(\log N)^\al$ and $\|\sum_{\la\in\La}\psi_\la\|\approx|\La|^{1/2}$.
\end{corollary}
\Proof We use the isomorphism $L^p\simeq L^p\oplus\ell^2$.
This isomorphism already implies that there is a conditional basis $\fX$ in $L^p$
with $k_N(\fX)\gtrsim N^\al$. Indeed, just take any unconditional basis $\sB_0$ in $L^p$
and a conditional basis $\sB_1$ in $\ell^2$ with $k_N(\sB_1)\gtrsim N^\al$ (using e.g. Proposition \ref{P1bis}), and consider in $L^p\oplus\ell^2$ the joint basis $\fX=\sB_0\cup\sB_1$
(say in alternate order).
This basis will satisfy $k_N(\fX)\gtrsim N^\al$ and also \eqref{xA}.

Now use that $L^p\simeq (L^p\oplus\ell^2)\oplus\ell^2$,
and apply the construction in Theorem \ref{ThO1}, using the basis $\fX$ in the first summand $\cX=L^p\oplus\ell^2$,
to obtain a new basis $\Psi$ with the required properties.
\ProofEnd

\section{Appendix I}

We show the claim asserted in Remark \ref{R1}.

\begin{proposition}
\label{P2}For $\al\in(0,\frac12)$, consider the Hilbert space
$L^2([-\pi,\pi],|t|^{-2\al}dt)$ with the trigonometric basis
$\{1,e^{it},e^{-it},e^{2it},e^{-2it},\ldots\}$. Then
\[
k_N\,\lesssim\,N^\al,\quad N=1,2,\ldots
\]
\end{proposition}

We shall use the following lemma.

\begin{lemma}
\label{Lvdc} Let $\ga\in(0,1)$ and $\om(t)=|t|^{-\ga}$, $|t|<\pi$.
Then its Fourier coefficients satisfy\Be |\widehat\om(k)|\, \leq\,
\frac{c_\ga}{|k|^{1-\ga}},\quad k\in\SZ\setminus\{0\}.
\label{wk}\Ee
\end{lemma}

\Proof We use the method of stationary phase. Let $\dt=1/|k|$.
Then
\[
\Big|\int_0^\dt e^{-ikt}\,t^{-\ga}\,dt\Big|\leq \int_0^\dt
t^{-\ga}\,dt =
\frac{\dt^{1-\ga}}{1-\ga}\,=\,\frac{c_\ga}{|k|^{1-\ga}}.
\]
On the other hand, by Van der Corput's Lemma \cite[p.
334]{Stein}\[ \Big|\int_\dt^\pi e^{-ikt}\,|t|^{-\ga}\,dt\Big|\leq
\frac{c}{|k|}\,\Big(\frac1{\pi^\ga}+\int_\dt^\pi \ga
|t|^{-\ga-1}\,dt\Big) =
\frac{c}{|k|\,\dt^\ga}\,=\,\frac{c}{|k|^{1-\ga}}.
\]
These two estimates easily imply \eqref{wk}.\ProofEnd

\

\Proofof{Proposition \ref{P2}} Let $|\La|=N$. It suffices to show
that \Be \big\|\sum_{\la\in \La} a_\la e^{i\la
t}\big\|_{L^2(|t|^{-2\al}dt)}^2\,\lesssim\,
N^{2\al}\,\sum_\la|a_\la|^2, \label{appaux1}\Ee since the basis is
besselian. We can write the left hand side as
\[
LHS =\int_{-\pi}^\pi\sum_{\la,\mu\in\La}a_\la\,\bar{a_\mu}\,
e^{i(\la-\mu)t}\,|t|^{-2\al}\,dt\,=\,\sum_{\la,\mu\in\La}a_\la\,\bar{a_\mu}\,\widehat\om(\mu-\la).
\]
Using Lemma \ref{Lvdc} (with $\ga=2\al$) we obtain \Bea LHS &
\lesssim &
\sum_\la|a_\la|^2\,+\,\sum_{\la\not=\mu\in\La}\frac{|a_\la|\,|a_\mu|}{|\mu-\la|^{1-2\al}}\nonumber\\
& \leq &
\sum_\la|a_\la|^2\,+\,\Big(\sum_\la|a_\la|^2\Big)^\frac12\,
\Big[\sum_{\la\in\La}\Big(\sum_{{\mu\in\La}\atop{\mu\not=\la}}\frac{|a_\mu|}{|\mu-\la|^{1-2\al}}\Big)^2\Big]^\frac12.\label{appaux2}\Eea
Now,\[
\sup_{\la\in\La}\,\sum_{{\mu\in\La}\atop{\mu\not=\la}}\frac{1}{|\mu-\la|^{1-2\al}}\,\leq\,
\sum_{j=1}^{|\La|}\frac1{j^{1-2\al}}\, \leq \, c N^{2\al},\] so
Schur's lemma (or Cauchy-Schwarz) easily leads from
\eqref{appaux2} to \eqref{appaux1}.\ProofEnd

\section{Appendix II}
 The following is a variation of Theorem \ref{Th1} for
$L^p$ spaces. Throughout this section we fix $1<p<\infty$, and let
$\|\cdot\|$ stand for the usual norm in $L^p(X,\mu)$. We denote by
$\kappa=\kappa(p)$ the smallest constant such that \Be
\|G_Nx\|\leq \kappa\|x\| \mand \|x-G_N x\|\leq \kappa\|x\|,
 \quad \forall\;x\in L^p,\;N=1,2,\ldots
\label{Kp}\Ee

\begin{theorem}
\label{Thp} Let $1<p<\infty$ and $\{\be_j\}$ a quasi-greedy
(normalized) basis in $L^p$. Then, there exists
$\al=\al(\kappa,p)<1$ and $c>0$ such that \Be
k_N\big(\{\be_j\}\big)\leq\, c\, (\log N)^\al,\quad
\forall\;N\in\SN. \label{kNalp}\Ee
\end{theorem}

The result depends on an $L^p$-version of Lemma \ref{L1}, which we
state with constants that very likely are not optimal.

\begin{lemma}\label{L1p}\

 \Benu \item[(i)] If $1<p\leq 2$
then for $c_p=2-\frac{p-1}{2\kappa^2}$ it holds
\Be\big\|x+y\big\|^2\,\leq\,c_p\,\Big(\|x\|^2+\|y\|^2\Big),\quad\forall\;x\succcurlyeq
y.\label{Par1}\Ee \item[(ii)] If $2\leq p<\infty$ then for
$c_p=2^{p-1}-\frac1{2\kappa^p}$ it holds
\Be\big\|x+y\big\|^p\,\leq\,c_p\,\Big(\|x\|^p+\|y\|^p\Big),\quad\forall\;x\succcurlyeq
y.\label{Par2}\Ee\Eenu
\end{lemma}

To prove Lemma \ref{L1p} we shall use weak versions of the
parallelogram identity which are well-known in the literature\Bea
\|x+y\|^2+(p-1)\|x-y\|^2\leq
2\big(\|x\|^2+\|y\|^2\big), & & 1<p\leq2\label{wpar1}\\
\|x+y\|^p+\|x-y\|^p\leq 2^{p-1}\big(\|x\|^p+\|y\|^p\big),& & 2\leq
p<\infty.\label{wpar2}\Eea The first one appears in work of Bynum
and Drew \cite{ByDr}, and the second one is attributed  to Clarkson \cite{Cla}.

\bline {\bf Case $1<p\leq 2$.}
 Call $N=|\supp x|$. The assumption $x\succcurlyeq y$ and
\eqref{Kp} give \Be \|x\|=\big\|G_N(x-y)\big\|\,\leq \,
\kappa\,\big\|x- y\big\|\mand \|y\|=\big\|(I-G_N)(x-
y)\big\|\,\leq \, \kappa\,\big\|x- y\big\|. \label{aux2}\Ee Thus,
\[ \big\|x- y\big\|^2\,\geq \,
\tfrac1{2\kappa^2}\,\big(\|x\|^2+\|y\|^2\big). \] Now the weak the
parallelogram law in \eqref{wpar1} combined with the previous
estimate gives \[ \|x+y\|^2\leq
2\big(\|x\|^2+\|y\|^2\big)\,-\,(p-1)\|x-y\|^2 \, \leq \,
\big(2-\tfrac{p-1}{2\kappa^2}\big)\,\big(\|x\|^2+\|y\|^2\big),\]
which proves \eqref{Par1}.

We now sketch the proof of Theorem \ref{Thp} in this case. Arguing
in exactly the same way as in the proof of Theorem \ref{Th1} one
reaches the inequality \eqref{aux1}, this time with $(1+\dt)$
replaced by the constant $c_p$. We follow the second approach
alluded in that proof, obtaining \eqref{aux1bis}, and hence the
validity of \eqref{aux0} with $\al=(1+\log_2 c_p)/2<1$. This
establishes \eqref{kNalp} when $1<p\leq 2$.

\bline {\bf Case $2\leq p<\infty$.} We first establish
\eqref{Par2}.  From \eqref{aux2} observe that \[\|x- y\|^p\,\geq
\, \tfrac1{2\kappa^p}(\|x\|^p+\|y\|^p).\] This, combined with
Clarkson's inequality in \eqref{wpar2} gives
\[ \|x+y\|^p\leq
2^{p-1}\big(\|x\|^p+\|y\|^p\big)\,-\,\|x-y\|^p \, \leq \,
\big(2^{p-1}-\tfrac{1}{2\kappa^p}\big)\,\big(\|x\|^p+\|y\|^p\big).\]
Now one proceeds as in the proof of Theorem \ref{Th1}, but with
powers 2 replaced by powers $p$, and $(1+\dt)$ replaced by $c_p$.
Then the corresponding version of \eqref{aux1bis} leads to\[
\|S_A(x)\|\lesssim\,c_p^{(\log_2 m)/{p}}\,m^{1/p}\,\|x\|\,=\,
m^{(1+\log_2 c_p)/p}\,\|x\|.
\]
So we can set $\al=(1+\log_2 c_p)/p$, which is smaller than 1
since $c_p<2^{p-1}$.

 \ProofEnd

 \

\bibliographystyle{plain}

\begin{thebibliography}{1}


\bibitem{Ba} \textsc {K.I. Babenko}, \emph{ On conjugate functions}. Dokl. Acad, Nauk SSSR 62, 157-160 (1948)
(in Russian).


\bibitem{ByDr} \textsc {W. L. Bynum and J. H. Drew}, \emph{A Weak Parallelogram Law for $\ell_p$},
Amer. Math. Monthly {\bf 79}(9) (1972), 1012--1015.

\bibitem{Cla} \textsc {J. A. Clarkson}, \emph{Uniformly Convex Spaces},
 Trans. Amer. Math. Soc. {\bf 40} (1936), 396-–414.

\bibitem{DKK}
\textsc{S.J. Dilworth, N.J. Kalton, D. Kutzarova},
\emph{On the existence of almost greedy bases in Banach spaces},  Studia Math. 159 (2003), no. 1, 67--101.

\bibitem{DKKT}
\textsc{S.J. Dilworth, N.J. Kalton, D. Kutzarova, and V.N.
Temlyakov}, \emph{The Thresholding Greedy Algorithm, Greedy Bases,
and Duality}, Constr. Approx., 19, (2003),575--597.

\bibitem{DST}
\textsc{S.J. Dilworth, M. Soto-Bajo, and V.N.
Temlyakov}, \emph{Quasi-greedy bases and Lebesgue-type
inequalities}, preprint 2012.

\bibitem{duo}
\textsc{J. Duoandikoetxea}, \emph{Fourier Analysis},
 Graduate Studies in Mathematics {\bf 29}, Amer. Math. Soc., Providence, 2001.

\bibitem{GHO}
\textsc{G. Garrig\'os, E. Hern\'andez and T. Oikhberg}, \emph{Lebesgue-type inequalities for quasi-greedy
bases}, preprint 2012.

\bibitem{Go}
\textsc{S. Gogyan}, \emph{An example of an almost greedy basis in
$L^1(0,1)$}, Proc. Amer. Math. Soc. {\bf 138} (4), (2010),
1425�-1432.

\bibitem{GG} \textsc{V.I. Gurariy, N.I. Gurariy}. \emph{Bases in uniformly convex and uniformly smooth Banach spaces}, Izv. Acad. Nauk. SSSR
ser. mat. {\bf 35} (1971) 210-215 (in Russian).

\bibitem{H}
     \textsc{E. Hern\'andez},
     \emph{Lebesgue-type inequalities for quasi-greedy bases}.
     Preprint 2011. ArXiv: 1111.0460v2 [matFA] 16 Nov 2011.

\bibitem{KT}
\textsc{S.V. Konyagin and V.N. Temlyakov}, \emph{A remark on greedy
approximation in Banach spaces}, East. J. Approx. 5, (1999),
365--379.

\bibitem{N1}
\textsc{M. Nielsen}, \emph{An example of an almost greedy
uniformly bounded orthonormal basis for $L^p(0, 1)$}, J. Approx.
Th. {\bf 149} (2007), 188--192.


\bibitem{Stein} \textsc{Stein, E.}, {\it Harmonic Analysis}. Princeton University Press, 1993.

\bibitem{Tem1}
\textsc{V.N. Temlyakov}, {\it Greedy approximation} Cambridge University Press


\bibitem{TYY1}
\textsc{V. N. Temlyakov, M. Yang, P. Ye}, \emph{Greedy approximation
with regard to non-greedy bases}, Adv. in Comp. Math., 34, (2011),
219--337.

\bibitem{TYY2}
\textsc{V. N. Temlyakov, M. Yang, P. Ye}, \emph{Lebesgue-type
inequalities for greedy approximation with respect to quasi-greedy
bases}, East J. Approx {\bf 17} (2011), 127--138.


\bibitem{Wo} \textsc{P. Wojtaszczyk},
\emph{Greedy Algorithm for General Biorthogonal Systems}, Journal of
Approximation Theory, 107, (2000), 293--314.

\end{thebibliography}

\vskip 1truemm

\end{document}